\documentclass[10pt,reqno]{amsart}
\usepackage{amsmath,amssymb,amsthm,amsfonts,amscd}
\usepackage{mathtools}
\usepackage{enumerate}
\usepackage{geometry}
\usepackage[colorlinks=true,linkcolor=blue,citecolor=red,urlcolor=blue]{hyperref}

\theoremstyle{plain}
\newtheorem{theorem}{Theorem}[section]
\newtheorem{lemma}[theorem]{Lemma}
\newtheorem{proposition}[theorem]{Proposition}
\newtheorem{corollary}[theorem]{Corollary}

\theoremstyle{definition}
\newtheorem{definition}[theorem]{Definition}

\theoremstyle{remark}
\newtheorem{remark}[theorem]{Remark}

\numberwithin{equation}{section}

\newcommand{\C}{\mathbb{C}}
\newcommand{\R}{\mathbb{R}}
\newcommand{\Z}{\mathbb{Z}}
\newcommand{\N}{\mathbb{N}}

\newcommand{\HH}{\mathbb{H}}          

\newcommand{\wt}{\operatorname{wt}}
\newcommand{\Aut}{\operatorname{Aut}}
\newcommand{\End}{\operatorname{End}}
\newcommand{\STr}{\operatorname{STr}}
\newcommand{\tr}{\operatorname{Tr}}
\newcommand{\SL}{\operatorname{SL}}
\newcommand{\Id}{\operatorname{Id}}

\newcommand{\T}{\mathcal{T}}        
\newcommand{\V}{\mathcal{V}}        

\newcommand{\bfx}{\mathbf{x}}
\newcommand{\bfv}{\mathbf{v}}
\newcommand{\bfz}{\mathbf{z}}
\newcommand{\bfy}{\mathbf{y}}

\renewcommand{\phi}{\varphi}
\newcommand{\vac}{\mathbf{1}_{V}}

\begin{document}
\title[Cohomology of Jacobi Forms]{Cohomology of Jacobi Forms}
\author{A. Zuevsky}
\address{Institute of Mathematics, Czech Academy of Sciences,
         \v{Z}itn\'{a} 25, Prague, Czech Republic}
\email{zuevsky@yahoo.com}

\begin{abstract}
We define and study a cohomology theory for the space of Jacobi $n$-point
functions generated by a vertex operator (super)algebra, using precise
analogues of Zhu's reduction formulas.
A cochain complex $(C^{\bullet}(W), \delta^{\bullet})$ is constructed whose
coboundary operators are given by Zhu-type reduction maps, and whose
cohomology groups $H^{n}_{J}(W)$ we call the  {reduction cohomology of
Jacobi forms}.
We prove that the $n$-th reduction cohomology of a $V$-module $W$ is
isomorphic to the space of analytic continuations of solutions to a
vertex-operator-algebraic analogue of the Knizhnik-Zamolodchikov equations.
We further show that Jacobi $n$-point reduction formulas are $n$-point
connections on the vertex operator algebra bundle over the torus, yielding
a Bott-Segal-type theorem: $H^{n}_{J}(W)$ is isomorphic to the cohomology
of the space of deformed sections of the VOA bundle.
\end{abstract}


\maketitle

\begin{center}
{Conflict of Interest and Data availability Statements:}
\end{center}

The author states that:

1.) The paper does not contain any potential conflicts of interests.

2.) The paper does not use any datasets. No datasets were generated
during and/or analysed during the current study.

3.) The paper includes all data generated or analysed during this study.

4.) Data sharing is not applicable to this article as no datasets were
generated or analysed during the current study.

5.) The data of the paper can be shared openly.

6.) No AI was used to write this paper.

\section{Introduction}
\label{SecIntro}

The computation of continuous cohomologies for non-commutative structures on
manifolds has proven to be a subject of significant geometrical interest
\cite{BS,Fei,Fuks,Wag}.
For Riemann surfaces, and even for higher-dimensional complex manifolds, the
classical cohomology of holomorphic vector fields is often trivial \cite{Kaw,Wag}.
Feigin \cite{Fei} obtained important results on (co)homology of cosimplicial
objects associated to holomorphic vector fields $\mathrm{Lie}(\mathfrak{M})$.

Vertex operator algebra (VOA) theory \cite{B,FHL,FLM,K} provides a rich source
of modular and quasi-modular forms through $n$-point functions (also called
$n$-point correlation or characteristic functions)
\cite{FHL,KZ,MT,MTZ,Zhu}.
These functions are subject to the action of differential operators with
specific analytic behaviour \cite{GK,GN,Ob} and satisfy functional equations
known as  {Zhu reduction formulas}, which express $(n+1)$-point functions
linearly in terms of $n$-point functions \cite{Zhu,MT,MTZ}.
Analogous formulas for Jacobi $n$-point functions (in which an element
$J\in V_{1}$ with semisimple $J(0)$-action is inserted to track a ``charge''
variable) were established by Bringmann-Krauel-Tuite \cite{BKT}.

\textbf{Aim and approach.}
In this paper we use the reduction formulas of \cite{BKT} as coboundary
operators in a cochain complex built from spaces of Jacobi $n$-point functions.
The resulting cohomology, which we call the  {reduction cohomology of
Jacobi forms}, simultaneously
encodes information about the non-commutative structure of $V$, 
produces functional equations for quasi-modular forms, and
admits a geometric interpretation in terms of connections on a VOA
        bundle over the torus.
Our approach is complementary to, and inspired by, Huang's work on first and
second cohomologies of grading-restricted vertex algebras \cite{Huang}, but it
operates at the level of $n$-point functions for all $n\geq 0$ simultaneously,
and incorporates the extra Jacobi (elliptic) variable.

\textbf{Plan and main results.}
Section \ref{chain} introduces the cochain complex
$(C^{\bullet}(W),\delta^{\bullet}(\cdot))$ and the chain condition
$\delta^{n+1}\circ\delta^{n}=0$, which amounts to a family of
functional-differential equations for the $n$-point functions
(Remark \ref{remnon-empty}).
Section \ref{SecZhu} gives explicit formulas for the coboundary operators in
four cases: the general case (Subsection \ref{main}), the simplest case
(Subsection \ref{simplest}), the shifted-Virasoro case
(Subsection \ref{shifted}), and the vertex operator superalgebra/orbifold case
(Subsection \ref{vosa}).
Section \ref{cohomology} contains the main results:
Theorem \ref{mainpropo} (the KZ characterisation of $H^n_J(W)$) and
Lemma \ref{pisaka} (the Bott-Segal-type theorem).
Appendices \ref{SecQuasi}, \ref{redya}, and \ref{durdom} recall, respectively,
the theory of (quasi-)Jacobi forms, the reduction formulas used as
coboundary operators, and the axioms of vertex operator (super)algebras.

The main results are:

\begin{theorem}[KZ characterisation - Theorem \ref{mainpropo}]
Under the assumptions of Subsections \ref{main}-\ref{vosa},
the $n$-th reduction cohomology $H^{n}_{J}(W)$
of a $V$-module $W$ is isomorphic, as a complex vector space, to
the space of analytic continuations of non-zero solutions
$\mathcal{Z}^{J}_{W}(\bfx_{n}; B)$ of the equation
\begin{equation}
\label{eqKZ_intro}
  \sum_{k=0}^{n}\sum_{m\geq 0}
  f_{k,m}(\bfx_{n};B) \;T_{k}(v_{n+1}[m]_{\beta} \;) \;
  \mathcal{Z}^{J}_{W}(\bfx_{n};B)=0,
  \quad x_{i}\notin\mathfrak{V}_{i},
\end{equation}
where $\beta=h$ for a shifted Virasoro element and $\beta=0$ otherwise.
These solutions are expressible as series of deformed Weierstrass functions
(Appendix \ref{defo}) recursively generated by the reduction formulas
\eqref{poros}.
\end{theorem}

\begin{theorem}[Bott-Segal theorem - Lemma \ref{pisaka}]
Jacobi $n$-point functions generated by the reduction formulas \eqref{poros}
are $n$-point connections on the space of $g$-automorphism-deformed sections of
the VOA bundle $\V$ over $\T$.
Consequently, $H^{n}_{J}(W)\cong{\mathcal Conn}^{n}/\mathcal{G}^{n-1}$ is isomorphic
to the cohomology of the space of deformed $\V$-sections.
\end{theorem}

The equation \eqref{eqKZ_intro} may be viewed as a Jacobi-form analogue of the
Knizhnik-Zamolodchikov (KZ) equations \cite{KZ,TK}.
The derivation of the KZ equations and that of the reduction formulas
both proceed via contour integration around auxiliary insertion points, followed
by application of the commutator formula for vertex operators (cf.\
Section \ref{SecZhu} and \cite[Sect. 4]{Y}).

Quasi-Jacobi forms appear in VOA theory in connection with
topological $\mathcal{N}=2$ vertex algebras \cite{HE},
Gromov-Witten potentials \cite{Kaw},
elliptic genera \cite{Lib},
and Landau-Ginzburg orbifolds \cite{KYY}.
The cohomology theory introduced here is a first step toward a complete
description of cohomologies of holomorphic objects arising from non-commutative
structures on complex manifolds.
It can be extended to higher-genus Riemann surfaces \cite{TZ} and applied to
the study of integrable models \cite{LS, RSZ} 
and deformation theory \cite{Ma}.

\section{Chain complex for Jacobi $n$-point functions}
\label{chain}
\subsection{Notation and Jacobi $n$-point functions}
\label{SSnotation}

Fix a vertex operator superalgebra $(V,Y,\vac,\omega)$ as recalled in
Appendix \ref{durdom}, and a $V$-module $W$.
For $n\geq 1$ and a tuple $\bfv_{n}=(v_{1},\ldots,v_{n})\in V^{\otimes n}$
we write $\bfv_{n}$ for both the tuple and the tensor product when context is
clear.
Mark $n$ points $p_{1},\ldots,p_{n}$ on the torus $\T=\C/(\Z+\tau\Z)$, with
local coordinates $z_{1},\ldots,z_{n}$ around them, and set
$\bfx_{n}=(\bfv_{n},\bfz_{n})$ with $\bfz_{n}=(z_{1},\ldots,z_{n})$.
Let $J\in V_{1}$ be such that $J(0)$ acts semisimply on $V$ (and on $W$).
Let $g\in\Aut(V)$ be an automorphism commuting with the parity automorphism
$\sigma$ defined by $\sigma a=(-1)^{p(a)}a$.

Setting $\widetilde{y}_{i}=(e^{z_{i}L(0)}v_{i},e^{z_{i}})$ and
$Y(\widetilde{\bfy}_{n})=Y(\widetilde{y}_{1})\cdots Y(\widetilde{y}_{n})$,
the  {orbifold Jacobi $n$-point function} associated to $W$ is
\begin{equation}
\label{defJacobi_npt}
  \mathcal{Z}^{J}_{W}(\bfx_{n};g,\tau)
  =\STr_{W}\;\bigl(Y(\widetilde{\bfy}_{n}) \;g \;q^{L(0)-c/24}\bigr),
  \quad q=e^{2\pi i\tau},
\end{equation}
where $\STr_{W}(X)=\tr_{W}(\sigma X)=\tr_{W_{\bar 0}}(X)-\tr_{W_{\bar 1}}(X)$.
The  {Jacobi $n$-point function} with charge parameter
$\zeta=e^{2\pi iz}$ is, for a weak $V$-module $W$ \cite{MT},
\begin{equation}
\label{defJacobi_npt_2}
  \mathcal{Z}^{J}_{W}(\bfx_{n};B)
  =\tr_{W}\;\bigl(Y(\widetilde{\bfy}_{n}) \;\zeta^{J(0)} \;q^{L(0)}\bigr),
\end{equation}
where $B$ denotes the parameters $(\tau,\zeta)$ (equivalently, $(\tau,z)$).
The one-point function reduces to
$\mathcal{Z}^{J}_{W}(x_{1};B)=\tr_{W}\;\bigl(o_{0}(v_{1}) \;\zeta^{J(0)} \;q^{L(0)}\bigr)$,
where $o_{0}(v)=v(\wt v-1)$ for $\wt v\in\Z$, and is independent of $z_{1}$
\cite{Zhu,BKT}.
The zero-point function is
$\mathcal{Z}^{J}_{W}(B)=\tr_{W}\;\bigl(\zeta^{J(0)} \;q^{L(0)}\bigr)$.
\subsection{The cochain complex}
\label{SScomplex}

For each $n\geq 0$ define the  {$n$-th cochain space}
\[
  C^{n}(W)=\bigl\{\mathcal{Z}^{J}_{W}(\bfx_{n};B)\mid
             \bfx_{n}\in V^{\otimes n}\times\C^{n},\;B=(\tau,z)\bigr\}.
\]
(With the convention that $C^{0}(W)=\{\mathcal{Z}^{J}_{W}(B)\}$.)

\begin{definition}[Coboundary operator]
\label{defcobdry}
For $n\geq 0$ and $\bfx_{n+1}=(v_{n+1},z_{n+1},\bfx_{n})$,
the  {coboundary operator}
$\delta^{n}(\bfx_{n+1}):C^{n}(W)\to C^{n+1}(W)$
is defined by
\begin{equation}
\label{poros}
  \delta^{n}(\bfx_{n+1}) \;\mathcal{Z}^{J}_{W}(\bfx_{n};B)
  =\sum_{k=0}^{n}\sum_{m\geq 0}
  f_{k,m}(\bfx_{n+1};B) \;T_{k}(v_{n+1}[m]_{\beta} \;) \;
  \mathcal{Z}^{J}_{W}(\bfx_{n};B),
\end{equation}
where $f_{k,m}(\bfx_{n+1};B)$ are meromorphic functions on $\T$ specified, 
in each of the four cases of Section \ref{SecZhu}, by the 
corresponding reduction formula recalled in Appendix \ref{redya}, 
$T_{k}(v_{n+1}[m]_{\beta} \;) \;\mathcal{Z}^{J}_{W}(\bfx_{n};B)
=\mathcal{Z}^{J}_{W}(T_{k}(v_{n+1}[m]_{\beta})\cdot\bfx_{n};B)$ 
is the insertion of the square-bracket mode $v_{n+1}[m]_{\beta}$
into the $k$-th slot of $\bfx_{n}$
        \[
          T_{k}(\Gamma)\cdot(x_{1},\ldots,x_{n})
          =(x_{1},\ldots,\Gamma\cdot x_{k},\ldots,x_{n});
        \]
$\beta=h$ (the shifted Virasoro element, Appendix \ref{squa}) in 
Subsection \ref{shifted}, and $\beta=0$ otherwise; and 
the sums in \eqref{poros} are finite for each fixed $n$-point function, 
as a consequence of the truncation property of vertex operators.
\end{definition}

\begin{remark}[Torsor interpretation]
The data $\bfx_{n}$ constitute a torsor for the group of transformations of
$V^{\otimes n}$ and of the local coordinates, in the sense of
\cite[Ch. 6]{BZF}.
Specifically, the operators $T_{k}(v_{n+1}[m]_{\beta})$ act on the
$V^{\otimes n}$-entries of $\bfx_{n}$, while the coefficient functions
$f_{k,m}(\bfx_{n+1};B)$ act on the coordinate part $\bfz_{n}$.
\end{remark}
\subsection{The chain condition and non-emptiness of cochains}

\begin{definition}[Admissible locus $\mathfrak{V}_{n}$]
For $n\geq 0$, let $\mathfrak{V}_{n}$ denote the subset of all $\bfx_{n}$
such that the  {chain condition}
\begin{equation}
\label{chain_cond}
  \delta^{n+1}(\bfx_{n+2}) \;\delta^{n}(\bfx_{n+1}) \;
  \mathcal{Z}^{J}_{W}(\bfx_{n};B)=0
\end{equation}
holds for all $\bfx_{n+2}$, $\bfx_{n+1}$, and all $B$.
\end{definition}

Expanding \eqref{chain_cond} using \eqref{poros} twice gives:

\begin{proposition}[Chain condition as functional equations]
\label{propchain}
The condition \eqref{chain_cond} is equivalent to the system
\begin{equation}
\label{conditions}
  \Biggl(\;\sum_{\substack{k'=0, \;k=0\\m',m\geq 0}}^{n+1, \;n}
  f_{k',m'}(\bfx_{n+2};B) \;f_{k,m}(\bfx_{n+1};B) \;
  T_{k'}(v_{n+2}[m']_{\beta}) \;T_{k}(v_{n+1}[m]_{\beta})\Biggr)
  \mathcal{Z}^{J}_{W}(\bfx_{n};B)=0
\end{equation}
holding for all choices of $\bfx_{n+2}$, $\bfx_{n+1}$.
\end{proposition}

\begin{remark}
\label{remnon-empty}
The relation \eqref{conditions} is an infinite family (indexed by $n\geq 0$)
of functional-differential equations, with finitely many terms for each fixed
$n$, on the holomorphic functions $\mathcal{Z}^{J}_{W}(\bfx_{n};B)$.
The coefficients $f_{k,m}(\bfx_{n+1};B)$ are generalised elliptic functions
as specified in Section \ref{SecZhu}.
The VOA elements $v_{k}\in V$ appear only through matrix elements and supertrace
insertions, so they do not appear explicitly in the final form of the equations.
By the general theory of such holomorphic functional-differential equations
\cite{FK,Gu}, each equation in the family \eqref{conditions} has a solution
in its domain of definition.
Hence the spaces $C^{n}(W)$ satisfying \eqref{conditions} are non-empty for
all $n\geq 0$.
\end{remark}

\begin{remark}
The relation \eqref{conditions} admits an analogue of Fay's trisecant identity
\cite{Fay} for vertex operator superalgebras, and can be interpreted as
producing new identities among quasi-modular forms on $\T$.
\end{remark}
\subsection{The reduction cohomology complex}

Restricting each $C^{n}(W)$ to the admissible locus $\mathfrak{V}_{n}$ (so
that $\delta^{n}\circ\delta^{n-1}=0$), one obtains a cochain complex
\begin{equation}
\label{complex}
  0 \longrightarrow C^{0}(W)
  \xrightarrow{\;\delta^{0}(x_{1})\;}
  C^{1}(W)
  \xrightarrow{\;\delta^{1}(\bfx_{2})\;}
  \cdots
  \xrightarrow{\;\delta^{n-1}(\bfx_{n})\;}
  C^{n}(W)
  \xrightarrow{\;\delta^{n}(\bfx_{n+1})\;}
  \cdots
\end{equation}

\begin{definition}[Reduction cohomology]
For $n\geq 1$, the  {$n$-th reduction cohomology of Jacobi forms} for
the $V$-module $W$ is
$H^{n}_{J}(W)
  ={\ker\delta^{n}(\bfx_{n+1})}/{\operatorname{Im}\delta^{n-1}(\bfx_{n})}$. 
\end{definition}

\section{Reduction formulas and coboundary operators}
\label{SecZhu}

In this section we make Definition \ref{defcobdry} explicit in four cases,
using the reduction formulas recalled in Appendix \ref{redya}.
Throughout, $J\in V_{1}$ and $J(0)v_{n+1}=\alpha v_{n+1}$ with $\alpha\in\C$.
\subsection{General coboundary operator}
\label{main}

We use Propositions \ref{propapnpt} and \ref{propapnpt0} from
Appendix \ref{redya}.
Assume $v_{n+1}\in V$ satisfies
$v_{n+1}[l]\cdot v_{k}=0$ for all $l\geq 1$ and $1\leq k\leq n$,
and $J(0)v_{n+1}=\alpha v_{n+1}$ with $\alpha\in\C$.
Summing Proposition \ref{propapnpt} (respectively \ref{propapnpt0})
over $l\in\Z$, multiplied by $z_{n+1}^{l-1}$, and using VOA associativity
(commutator formula \eqref{eqcomm}), the coboundary operator \eqref{poros}
takes the form
\begin{equation}
\label{eqgeneral_cobdry}
  \delta^{n}(\bfx_{n+1}) \;\mathcal{Z}^{J}_{W}(\bfx_{n};z,\tau)
  =\sum_{\substack{l\in\Z, \;m\geq 0 \\ k=0}}^{n}
  f_{k,m}(\bfx_{n+1};B) \;T_{k}(v_{n+1}[m]_{\beta}) \;
  \mathcal{Z}^{J}_{W}(\bfx_{n};z,\tau),
\end{equation}
with coefficient functions
\begin{align}
  f_{0}(\bfx_{n+1};B) \;T_{0}(v[m])
  &=\sum_{l\in\Z}(-1)^{l+1} \;\delta_{\alpha z, \;\Z\tau+\Z} \;
    \frac{\lambda^{l-1}}{(l-1)!} \;z_{n+1}^{l-1} \;T_{0}(o_{\lambda}(v_{n+1})),
  \label{eqf0_general}\\[4pt]
  f_{k,m}(\bfx_{n+1};B)
  &=\sum_{l\in\Z}(-1)^{m+1}\binom{m+l-1}{m}z_{n+1}^{l-1}
    F_{k,m}(\bfx_{n+1};l,\alpha z,\tau),
  \label{persina}
\end{align}
where $\delta_{\alpha z,\Z\tau+\Z}=1$ if $\alpha z\in\Z\tau+\Z$ and $0$
otherwise, and
\begin{align*}
  F_{k,m}(\bfx_{n+1};l,\alpha z,\tau)
  &=\delta_{0,m} \;T^{1-\delta_{\alpha z,\Z\tau+\Z}}
    \;. \;\widetilde{E}_{m+l,\lambda}
    \;\bigl((1-\delta_{\alpha z,\Z\tau+\Z})\alpha z,\tau\bigr)\\
  &\quad +T^{1-\delta_{\alpha z,\Z\tau+\Z}}
    \;. \;\widetilde{P}_{m+l,(1-\delta_{\alpha z,\Z\tau+\Z})\lambda}
    \;\Bigl(\frac{z_{1}-z_{k}}{2\pi i},
             (1-\delta_{\alpha z,\Z\tau+\Z})\alpha z,\tau\Bigr),
\end{align*}
with $\widetilde{E}_{m+k,\lambda}(\alpha z,\tau)$ and
$\widetilde{P}_{m+l,\lambda}(z',\alpha z,\tau)$ as given in \eqref{eqGkl}
and \eqref{eqPellPm}, respectively.
Here the tilde denotes the action of the operator $T$ on the corresponding
special function: $T.E_{m+l,\lambda}=\widetilde{E}_{m+l,\lambda}$, etc.
\subsection{Simplest coboundary operator}
\label{simplest}

Suppose $J(0)v_{n+1}=\alpha v_{n+1}$ with $\alpha\in\C$.
Using Propositions \ref{propZhured} and \ref{propZhured0} directly (without
summing over $l\in\Z$), the coboundary operator \eqref{poros} is given by
\begin{eqnarray*}
  f_{0}(\bfx_{n+1};\alpha z,\tau) \;T_{0}(v_{n+1}[m])
  =\delta_{\alpha z, \;\lambda\tau+\mu\in\Z\tau+\Z} \;
    e^{-z_{n+1}\lambda} \;T_{0}(o_{\lambda}(v_{n+1})),
\label{eqf0_simple}
\qquad \qquad \qquad &&
\\
  f_{k,m}(z_{n+1};\lambda,k,\alpha z,\tau)
  =T^{1-\delta_{\alpha z, \;\lambda\tau+\mu\in\Z\tau+\Z}}\;
    P_{m+1,\lambda}
    \;\Bigl(\frac{z_{n+1}-z_{k}}{2\pi i},
             (1-\delta_{\alpha z, \;\lambda\tau+\mu\in\Z\tau+\Z})\alpha z,\tau\Bigr),&&
  \label{eqfk_simple}
\end{eqnarray*}
with $P_{m+1,\lambda}(w,\tau)$ defined in \eqref{eqPellPm}.
\subsection{Coboundary operator for a shifted Virasoro vector}
\label{shifted}

Suppose $J(0)a=\alpha a$ with $\alpha\notin\Z\setminus\{0\}$, and let
$g=e^{2\pi i(\mu/\alpha)J(0)}\in\Aut(V)$ for $\mu\in\Z$ chosen so that
$ga=a$.
Setting $J(0)v_{k}=\alpha_{k}v_{k}$ for $k=1,\ldots,n$, and using the shifted
Virasoro grading $L_{h}(0)=L(0)+(\lambda/\alpha)J(0)$ (Appendix \ref{squa}),
define $\widetilde{\widetilde{y}}_{i}=(e^{z_{i}L_{h}(0)}v_{i},e^{z_{i}})$.
The  {shifted Jacobi form}
\[
  \mathcal{Z}^{J}_{W}(\bfx_{n+1};h,\mu,\alpha,z,\tau)
  =\tr_{W}\;\bigl(Y(\widetilde{\widetilde{\bfy}}_{n+1}) \;g \;q^{L_{h}(0)}\bigr)
\]
is annihilated by $\sum_{k=1}^{n}\tr_{W}(T_{k}(a[0])Y(\widetilde{\bfy}_{n}) \;g \;q^{L(0)})=0$ (Corollary \ref{corZeroRes}).
The coboundary operator is \eqref{poros} with
\begin{eqnarray*}
  f_{0}(\bfx_{n+1};B) \;T_{0}(v_{n+1}[m]_h)=T_{0}(o_{h}(v_{n+1})),
  \\
  f_{k,m}(\bfx_{n+1};B)=P_{m+1}\;\Bigl(\frac{z_{n+1}-z_{k}}{2\pi i},\tau\Bigr),
\end{eqnarray*}
where $o_{h}(v_{n+1})=v_{n+1}(\wt_{h}(v_{n+1})-1)
=v_{n+1}(\wt(v_{n+1})-1+\mu)=o_{\mu}(v_{n+1})$.
\subsection{Vertex operator superalgebra (orbifold) case}
\label{vosa}

Let $v_{n+1}$ be homogeneous of weight $\wt(v_{n+1})\in\R$, and set
$\phi=\exp(2\pi i \;\wt(v_{n+1}))\in U(1)$.
Suppose $v_{n+1}$ is a $g$-eigenstate with $gv_{n+1}=\theta^{-1}v_{n+1}$
for some $\theta\in U(1)$, so that $g^{-1}v_{n+1}(k)g=\theta v_{n+1}(k)$.
The coboundary operator is given by \eqref{poros} with
\begin{eqnarray*}
  f_{0}(\bfx_{n+1};B) \;T_{0}(v_{n+1}[m])
  =\delta_{\theta,1} \;\delta_{\phi,1} \;T_{0}(o_{0}(v_{n+1})),
  \\
  f_{k,m}(\bfx_{n+1};B)
  =p(v_{n+1},\bfv_{k-1}) \;P_{m+1}\;\begin{bmatrix}\theta\\\phi\end{bmatrix}\;
   (z_{n+1}-z_{k},\tau),
\end{eqnarray*}
where the deformed Weierstrass function $P_{m+1}\;\bigl[\begin{smallmatrix}
\theta\\\phi\end{smallmatrix}\bigr]$ is defined in \eqref{Pkuv} and
$p(v_{n+1},\bfv_{k-1})$ is the sign factor of Proposition \ref{Propa[0]comm}.
Note that this case is related to the shifted-Virasoro case
(Subsection \ref{shifted}) upon specialising $\theta=\phi=1$.

\section{Cohomology}
\label{cohomology}
\subsection{The \texorpdfstring{$n$}{n}-th cohomology and the KZ equations}

We now state and prove the main theorem.

\begin{theorem}[KZ characterisation of $H^{n}_{J}(W)$]
\label{mainpropo}
Under the notations and assumptions of Subsections \ref{main}-\ref{vosa},
the $n$-th reduction cohomology $H^{n}_{J}(W)$ of the $V$-module $W$
is isomorphic, as a $\C$-vector space, to the space of analytic continuations
of non-zero solutions $\mathcal{Z}^{J}_{W}(\bfx_{n};B)$
to the equation
\begin{equation}
\label{poroserieroj}
  \sum_{k=0}^{n}\sum_{m\geq 0}
  f_{k,m}(\bfx_{n};B) \;T_{k}(v_{n+1}[m]_{\beta}) \;
  \mathcal{Z}^{J}_{W}(\bfx_{n};B)=0,
  \qquad x_{i}\notin\mathfrak{V}_{i},\;1\leq i\leq n,
\end{equation}
where $\beta=h$ in the shifted-Virasoro case and $\beta=0$ otherwise.
These solutions are quasi-modular forms expressible as series in deformed
Weierstrass functions (Appendix \ref{defo}), recursively generated by the
reduction formulas \eqref{poros}; their analytic continuations extend outside
$\mathfrak{V}_{n}$.
\end{theorem}

\begin{remark}
Equation \eqref{poroserieroj} is a vertex-operator-algebraic analogue of the
Knizhnik-Zamolodchikov equations \cite{KZ,TK} in the setting of Jacobi forms.
\end{remark}

\begin{remark}
Theorem \ref{mainpropo} relates the cohomological structure of Jacobi forms,
defined through VOA $n$-point functions, to the analytic structure of solutions
of a fundamental equation in mathematical physics.
For the cases $n=1,2$ one recovers, via Huang's results \cite{Huang}, the
first and second cohomologies of a grading-restricted vertex algebra in terms
of derivations and square-zero extensions.
\end{remark}

\begin{proof}
\textit{Representatives of $\ker\delta^{n}$.} 
The kernel of $\delta^{n}(\bfx_{n+1})$ acting on $C^{n}(W)$
consists precisely of those $\mathcal{Z}^{J}_{W}(\bfx_{n};B)$
satisfying \eqref{poroserieroj}, for $\bfx_{n}$ with $x_{i}\notin\mathfrak{V}_{i}$.
(For $x_{i}\in\mathfrak{V}_{i}$ the function vanishes by definition of the
admissible locus.)

\textit{Image of $\delta^{n-1}$ and recursion.}
The image $\operatorname{Im}\delta^{n-1}(\bfx_{n})$ consists of $n$-point functions
$\mathcal{Z}^{J}_{W}(\bfx'_{n};B)$ of the form
\begin{equation}
\label{poroserieroj_2}
  \mathcal{Z}^{J}_{W}(\bfx'_{n};B)
  =\Biggl(\sum_{k=1}^{n-1}\sum_{m\geq 0}
   f_{k,m}(\bfx_{n};B) \;T^{(g)}_{k}(v'_{n}[m]_{\beta})\Biggr)
   \mathcal{Z}^{J}_{W}(\bfx'_{n-1};B).
\end{equation}
Applying the reduction formulas \eqref{poros} recursively for each $x_{i}$,
$1\leq i\leq n$, with $x_{i}\notin\mathfrak{V}_{i}$, one obtains
\begin{equation}
\label{topaz}
  \mathcal{Z}^{J}_{W}(\bfx_{n};B)
  =\mathcal{D}(\bfx_{n};B) \;\mathcal{Z}^{J}_{W}(B),
\end{equation}
where $\mathcal{D}(\bfx_{n};B)$ is a differential operator built from the
coefficients $f_{k,m}$ and the modes $v_{i}[m]$; see \cite{MT,TZ} for explicit
formulas in the torus case.
If $x_{i}\in\mathfrak{V}_{i}$ at any step, the recursion forces
$\mathcal{Z}^{J}_{W}(\bfx_{n};B)=0$.
Analogously, $\mathcal{Z}^{J}_{W}(\bfx'_{n};B)$ given by \eqref{poroserieroj_2}
either vanishes (when $v_{n-i}\in\mathfrak{V}_{n-i}$ for some
$2\leq i\leq n$) or equals \eqref{topaz} with arguments $\bfx'_{n}$.

\textit{KZ interpretation.}
The derivation of the reduction relations \eqref{poros} in \cite{Y,BKT}
proceeds by double contour integration of
$\mathcal{Z}^{J}_{W}(\bfx_{n};B)$ around auxiliary variables, with
appropriate reproduction kernels, in exactly the same way that the
KZ equations are derived from the conformal Ward identities \cite{KZ,TK}.
Thus \eqref{poroserieroj} is indeed a KZ-type equation; explicit solutions for
specific $V$ and torus configurations appear in \cite{Y,MT}.

\textit{Analytic continuation.}
Using the commutator formula \eqref{eqcomm} and VOA associativity, the action
of $T_{k}(v_{n+1}[m]_{\beta})$ on $v_{k}\in V$ can be transferred to a shift
$z_{i}\mapsto z'_{i}=z_{i}+z_{n+1}$, $1\leq i\leq n$, of the coordinate
parameters, without changing the $V$-part $\bfv_{n}$ of $\bfx_{n}$.
Convergence of $n$-point functions in the shifted domain $z_{n+1}+\mathcal{T}_{n}$
(where $\mathcal{T}_{n}\subset\T$ is the locus $z_{i}\neq z_{j}$ for $i\neq j$)
follows from standard estimates on torus $n$-point functions \cite{Zhu,MT}.
Hence elements of $H^{n}_{J}(W)$ are analytic continuations of solutions to
\eqref{poroserieroj} across $\mathfrak{V}_{n}$.
\end{proof}

\begin{corollary}
\label{cororbifold_det}
In the orbifold case (Subsection \ref{vosa}) with $v_{n}\in V_{n}$,
the $n$-th cohomology $H^{n}_{J}(W)$ contains the space of functions
of the form $\det\bigl(M(z_{i}-z_{j})_{1\leq i,j\leq n}\bigr)$,
where $M(z_{i}-z_{j})$ is an $n\times n$ matrix whose entries are
deformed elliptic functions depending on $z_{i}-z_{j}$,
running over all possible mode combinations of $v_{n}$.
\end{corollary}

\begin{proof}
Apply Theorem \ref{Theorem_npt_rec0} iteratively, using the fact that the
deformed Weierstrass functions $P_{m+1}\bigl[\begin{smallmatrix}\theta\\
\phi\end{smallmatrix}\bigr](z_{i}-z_{j},\tau)$ combine under the
anti-symmetrisation implied by the supertrace to produce determinants
of the stated form \cite{MTZ,BKT}.
\end{proof}
\subsection{Geometric meaning: connections and a Bott-Segal theorem}
\label{meme}

We now show that the reduction formulas are $n$-point connections
on a VOA bundle, and derive a Bott-Segal-type theorem for $H^{n}_{J}(W)$.

\begin{definition}[Multi-point connection]
\label{defnpt_conn}
Let $\V$ be a holomorphic vector bundle over $\T$, and let $\mathcal{T}_{0}
\subset\T$ be an open subdomain.
Denote by $\mathcal{SV}$ the sheaf of holomorphic sections of $\V$.
A  {multi-point connection} on $\V$ is a $\C$-multilinear map
$\mathcal{G}$ such that, for any holomorphic function $f$ on $\mathcal{T}_{0}$
and sections $\phi\in\mathcal{SV}(p)$, $\psi\in\mathcal{SV}(p')$ at points
$p,p'\in\mathcal{T}_{0}$,
\begin{equation}
\label{locus}
  \sum_{q,q'\in\mathcal{T}_{0}}
  \mathcal{G}\bigl(f(\psi(q))\cdot\phi(q')\bigr)
  =f(\psi(p')) \;\mathcal{G}(\phi(p))
  +f(\phi(p)) \;\mathcal{G}(\psi(p')),
\end{equation}
where the left-hand sum converges absolutely on $\mathcal{T}_{0}$.
The space of all such $n$-point connections is denoted ${\mathcal Conn}^{n}$.
\end{definition}

The  {form} of a multi-point connection $\mathcal{G}$ is
\begin{equation}
\label{gform}
  G(f,\phi,\psi)
  =f(\phi(p)) \;\mathcal{G}(\psi(p'))
  +f(\psi(p')) \;\mathcal{G}(\phi(p))
  -\sum_{q,q'\in\mathcal{T}_{0}}
   \mathcal{G}\bigl(f(\psi(q'))\cdot\phi(q)\bigr),
\end{equation}
and the space of connection forms is $\mathcal{G}^{n}$.

\begin{lemma}[Bott-Segal theorem for Jacobi forms]
\label{pisaka}
The Jacobi $n$-point forms \eqref{defJacobi_npt_2}, generated by the
reduction formulas \eqref{poros}, are $n$-point connections on the space of
$g$-automorphism-deformed sections of the VOA bundle $\V$ over $\T$
(the Virasoro-shifted version of the bundle constructed in \cite[Ch. 6]{BZF}).
For $n\geq 0$, the $n$-th reduction cohomology satisfies
\[
  H^{n}_{J}(W)\cong{\mathcal Conn}^{n}/\mathcal{G}^{n-1},
\]
i.e., it is isomorphic to the cohomology of the space of deformed $\V$-sections.
\end{lemma}

\begin{remark}
Lemma \ref{pisaka} is the analogue, for deformed VOA-bundle sections on the
torus, of the main result of Bott-Segal \cite{BS} and Wagemann \cite{Wag}
on holomorphic vector fields.
\end{remark}

\begin{proof}
\textit{Identification of $n$-point functions with connections.}
The VOA bundle $\V$ and its $g$-twisted dual $\V^{*}$ were explicitly
constructed in \cite[Ch. 6, \S6.5.3]{BZF}.
We use the Virasoro-shifted version, replacing $L(0)$ by $L_{h}(0)$ throughout.
The intrinsic (coordinate-independent) VOA operators are defined by
\[
  \langle u, \;(\mathcal{Y}^{*}_{\mathbf{p}}(i(\bfv_{n})))_{n} \;g \;v\rangle
  =\langle u, \;Y(\bfx_{n}) \;v\rangle,
\]
realised as matrix elements of vertex operators on punctured disks with local
coordinates $\bfz_{n}$ on $\T$ (cf.\ \cite[Prop. 6.5.4]{BZF}).

\textit{Setting up the identification.}
For non-vanishing $f(\phi(p))$, define the identification
\begin{align*}
  \mathcal{G}                         &= \mathcal{Z}^{J}_{W}(\bfx_{n};B),\\
  \psi(p')                            &= \bfx_{n+1},\quad
  \phi(p) = \bfx_{n},\\
  \mathcal{G}(f(\psi(q))\cdot\phi(q')) &= T_{k}(v[m]_{\beta}) \;\mathcal{Z}^{J}_{W}(\bfx_{n};B),\\
  -\tfrac{f(\psi(p'))}{f(\phi(p))} \;\mathcal{G}(\phi(p))
                                      &= f_{0}(\bfx_{n+1};B) \;T_{0}(o_{\lambda}(v_{n+1})) \;
                                         \mathcal{Z}^{J}_{W}(\bfx_{n};B),\\
  f^{-1}(\phi(p))\;
  \sum_{q_{n},q'_{n}\in\mathcal{T}_{0}}
  \mathcal{G}(f(\psi(q))\cdot\phi(q'))
                                      &= \sum_{k=1}^{n}\sum_{m\geq 0}
                                         f_{k,m}(\bfx_{n+1};B) \;
                                         T_{k}(v[m]_{\beta}) \;
                                         \mathcal{Z}^{J}_{W}(\bfx_{n};B).
\end{align*}
Under these identifications, the multi-point connection property \eqref{locus}
becomes exactly the reduction formula \eqref{poros}.

\textit{Cohomology identification.}
By \cite[Constr. 6.6.4 and Prop. 6.6.7]{BZF}, the $n$-point functions are
holomorphic connections on $\V$, and the reduction cohomology quotient
$\ker\delta^{n}/\operatorname{Im}\delta^{n-1}$ coincides with
${\mathcal Conn}^{n}/\mathcal{G}^{n-1}$, the cohomology of the space of
deformed $\V$-sections.
\end{proof}
\subsection{Geometric meaning of the chain condition \eqref{conditions}}

Since the operators in \eqref{poros} act only on the $V$-part of $\bfx_{n}$,
the condition \eqref{conditions} has the following geometric interpretation.
All operators $T_{k}$ change VOA elements via either the zero-mode
$o(v)=v_{\wt v-1}$ or positive square-bracket modes $v[m]$, $m\geq 0$.
Since $n$-point Jacobi forms are quasi-modular forms \cite{BKT},
the condition \eqref{conditions} encodes new relations among quasi-modular
forms; it also defines a complex-analytic subvariety in $\bfz_{n}$
with non-commutative parameters $\bfv_{n}\in V^{\otimes n}$.
The modular-invariance of higher-$n$ Jacobi forms can be proved using
\eqref{poros}, analogously to \cite{MT,MTZ}.

\appendix
\section{Quasi-Jacobi forms}
\label{SecQuasi}
\subsection{Jacobi forms}
\label{SSJac}

We recall definitions from \cite{EZ,BKT}.
Let $k,m\in\N_{0}$ and let $\chi$ be a rational character of the Jacobi group
$\SL(2,\Z)\ltimes\Z^{2}$.
A  {holomorphic Jacobi form} of weight $k$ and index $m$ on $\SL_{2}(\Z)$
with multiplier $\chi$ is a holomorphic function $\phi:\C\times\HH\to\C$
satisfying: for $\gamma=\bigl(\begin{smallmatrix}a&b\\c&d\end{smallmatrix}\bigr)
\in\SL_{2}(\Z)$ and $(\lambda,\mu)\in\Z\times\Z$,
\begin{equation}
\label{eqJac_transform}
  \phi\Bigl|_{k,m}\bigl(\gamma,(\lambda,\mu)\bigr)=\chi\bigl(\gamma,(\lambda,\mu)\bigr)\phi,
\end{equation}
where
\begin{eqnarray*}
  \phi\Bigl|_{k,m}\bigl(\gamma,(\lambda,\mu)\bigr)(z,\tau)
  =\frac{     
e\left(-\frac{cm(z+\lambda\tau+\mu)^{2}}{c\tau+d}
   +m(\lambda^{2}\tau+2\lambda z)\right)}{(c\tau+d)^{k}}\;
   \phi\left(\frac{z+\lambda\tau+\mu}{c\tau+d},\gamma.\tau\right),
\end{eqnarray*}
with $e(w)=e^{2\pi iw}$.
Writing $\chi\bigl(\begin{smallmatrix}1&1\\0&1\end{smallmatrix}\bigr)
=e^{2\pi ia_{1}/N_{1}}$ and $\chi(0,1)=e^{2\pi ia_{2}/N_{2}}$ for uniquely
determined $N_{j}\in\N$, $a_{j}\in\N$ with $\gcd(a_{j},N_{j})=1$,
the Fourier expansion of $\phi$ is
\[
  \phi(z,\tau)
  =\sum_{\substack{n\in\N_{0}+\rho_{1}\\r\in\Z+\rho_{2},\;r^{2}\leq 4nm}}
   c(n,r) \;q^{n}\zeta^{r},
  \quad q=e({\tau}),\;\zeta=e(z),\;\rho_{j}=\tfrac{a_{j}}{N_{j}}\bmod\Z.
\]
\subsection{Quasi-Jacobi forms}

Following \cite{Lib}, an  {almost meromorphic Jacobi form} of weight $k$,
index $0$, and depth $(s,t)$ is a meromorphic function in
$\C\{q,\zeta\}[z^{-1},z_{2}/\tau_{2},1/\tau_{2}]$
(where $z=z_{1}+iz_{2}$, $\tau=\tau_{1}+i\tau_{2}$) satisfying the
transformation law \eqref{eqJac_transform} and of degree at most $s$ (resp.\
$t$) in $z_{2}/\tau_{2}$ (resp.\ $1/\tau_{2}$).
A  {quasi-Jacobi form} of weight $k$, index $0$, and depth $(s,t)$ is the
constant term (in $z_{2}/\tau_{2}$ and $1/\tau_{2}$) of such a form.
\subsection{Modular and elliptic functions}
\label{SSell}

For a variable $x$ set $D_{x}=\tfrac{1}{2\pi i}\tfrac{\partial}{\partial x}$
and $q_{x}=e^{2\pi ix}$.
For $m\in\N=\{\ell\in\Z:\ell>0\}$, the  {Weierstrass functions} are
\begin{gather*}
  P_{1}(w,\tau)=-\sum_{n\in\Z\setminus\{0\}}\frac{q_{w}^{n}}{1-q^{n}}-\tfrac{1}{2},\\
  P_{m+1}(w,\tau)=\frac{(-1)^{m}}{m!}D_{w}^{m}P_{1}(w,\tau)
  =\frac{(-1)^{m+1}}{m!}\sum_{n\in\Z\setminus\{0\}}\frac{n^{m}q_{w}^{n}}{1-q^{n}}.
\end{gather*}
The  {Eisenstein series} for $k\in 2\N$ are
\[
  E_{k}(\tau)=-\frac{B_{k}}{k!}+\frac{2}{(k-1)!}\sum_{n\geq 1}\frac{n^{k-1}q^{n}}{1-q^{n}},
\]
where $B_{k}$ is the $k$-th Bernoulli number ($(e^{z}-1)^{-1}=\sum_{k\geq 0}B_{k}z^{k-1}/k!$),
and $E_{k}=0$ for odd $k$.
We set $E_{0}=-1$.
Then $E_{k}$ is a modular form of weight $k$ for $k>2$, and quasi-modular for
$k=2$:
\[
  E_{k}(\gamma\tau)=(c\tau+d)^{k}E_{k}(\tau)-\delta_{k,2} \;\frac{c(c\tau+d)}{2\pi i}.
\]
For $w,z,\tau\in\HH$ the  {twisted Weierstrass functions} are
\begin{gather*}
  \widetilde{P}_{1}(w,z,\tau)=-\sum_{n\in\Z}\frac{q_{w}^{n}}{1-q_{z}q^{n}},\\
  \widetilde{P}_{m+1}(w,z,\tau)=\frac{(-1)^{m+1}}{m!}\sum_{n\in\Z}\frac{n^{m}q_{w}^{n}}{1-q_{z}q^{n}}.
\end{gather*}
For $m\in\N_{0}$ and $\lambda\in\C$, the  {$\lambda$-twisted Weierstrass
functions} are \cite{Zag}
\begin{equation}
\label{eqPellPm}
  P_{m+1,\lambda}(w,\tau)
  =\frac{(-1)^{m+1}}{m!}\sum_{n\in\Z\setminus\{-\lambda\}}\frac{n^{m}q_{w}^{n}}{1-q^{n+\lambda}}.
\end{equation}
One has $P_{1,\lambda}(w,\tau)=q_{w}^{-\lambda}(P_{1}(w,\tau)+\tfrac{1}{2})$
and $P_{m+1,\lambda}(w,\tau)=\frac{(-1)^{m}}{m!}D_{w}^{m}P_{1,\lambda}(w,\tau)$.
The Laurent expansion
\[
  P_{1,\lambda}(w,\tau)=\frac{1}{2\pi iw}-\sum_{k\geq 1}E_{k,\lambda}(\tau)(2\pi iw)^{k-1}
\]
defines the  {twisted Eisenstein series} \cite{Zag}
\begin{equation}
\label{eqGkl}
  E_{k,\lambda}(\tau)=\sum_{j=0}^{k}\frac{\lambda^{j}}{j!}E_{k-j}(\tau).
\end{equation}
The alternative generating function 
\[
\widetilde{P}_{1}(w,z,\tau)
=\frac{1}{2\pi iw}-\sum_{k\geq 1}\widetilde{E}_{k}(z,\tau)(2\pi iw)^{k-1}, 
\]
\cite{Ob} gives, for $k\geq 1$,
\[
  \widetilde{E}_{k}(z,\tau)
  =-\delta_{k,1}\frac{q_{z}}{q_{z}-1}-\frac{B_{k}}{k!}
  +\frac{1}{(k-1)!}\sum_{m,n\geq 1}
   \bigl(n^{k-1}q_{z}^{m}+(-1)^{k}n^{k-1}q_{z}^{-m}\bigr)q^{mn},
\]
and $\widetilde{E}_{0}(z,\tau)=-1$.
\subsection{Deformed elliptic functions}
\label{defo}

Following \cite{DLM,MTZ}, let $(\theta,\phi)\in U(1)\times U(1)$ with
$\phi=\exp(2\pi i\lambda)$, $0\leq\lambda<1$.
For $z\in\C$, $\tau\in\HH$, the  {deformed Weierstrass functions} are,
for $k\geq 1$,
\begin{equation}
\label{Pkuv}
  P_{k}\;\begin{bmatrix}\theta\\\phi\end{bmatrix}\;(z,\tau)
  =\frac{(-1)^{k}}{(k-1)!}
  \sideset{}{'}\sum_{n\in\Z+\lambda}
  \frac{n^{k-1}q_{z}^{n}}{1-\theta^{-1}q^{n}},
\end{equation}
where the prime means that the $n=0$ term is omitted if $(\theta,\phi)=(1,1)$.
These functions converge absolutely and uniformly on compact subsets of
$\{|q|<|q_{z}|<1\}$ \cite{DLM}, and satisfy
\[
P_{k}\bigl[\begin{smallmatrix}\theta\\\phi\end{smallmatrix}\bigr](z,\tau)
=\frac{(-1)^{k-1}}{(k-1)!}\frac{d^{k-1}}{dz^{k-1}}
 P_{1}\bigl[\begin{smallmatrix}\theta\\\phi\end{smallmatrix}\bigr](z,\tau).
\]
\section{Reduction formulas for Jacobi $n$-point functions}
\label{redya}
We recall the reduction formulas from \cite{BKT,MTZ}.
\subsection{Superalgebra case}

\begin{proposition}[{\cite[Prop. 3.1]{BKT}}]
\label{Propa[0]comm}
Let $v_{n+1}\in V$ be homogeneous of integer weight $\wt(v_{n+1})\in\Z$.
Then
\[
  \sum_{k=1}^{n}p(v_{n+1},\bfv_{k-1}) \;
  \mathcal{Z}^{J}_{W}((v_{n+1}[0]\cdot)_{k}\bfv_{n};B)=0,
\]
where $p(A,B_{1}\cdots B_{r-1})=(-1)^{p(A)(p(B_{1})+\cdots+p(B_{r-1}))}$
for $r>1$ and $p(A,\emptyset)=1$ for $r=1$.
\end{proposition}

Let $\phi=\exp(2\pi i \;\wt(v_{n+1}))\in U(1)$ and
$gv_{n+1}=\theta^{-1}v_{n+1}$ with $g^{-1}v_{n+1}(k)g=\theta v_{n+1}(k)$
for some $\theta\in U(1)$.

\begin{theorem}[{\cite[Thm. 3.2]{BKT}}]
\label{Theorem_npt_rec0}
Let $v_{n+1}$, $\theta$, $\phi$ be as above.
For any $\bfv_{n}\in V^{\otimes n}$,
\begin{align*}
  \mathcal{Z}^{J}_{W}(\bfx_{n+1};B)
  &=\delta_{\theta,1} \;\delta_{\phi,1} \;
    \STr_{W}\;\bigl(o(v_{n+1}) \;Y_{W}(\bfy_{n}) \;g \;q^{L(0)-c/24}\bigr)\\
  &+\sum_{k=1}^{n}\sum_{m\geq 0}p(v_{n+1},\bfv_{k-1}) \;
    P_{m+1}\;\begin{bmatrix}\theta\\\phi\end{bmatrix}\;
    (z_{n+1}-z_{k},\tau) \;
    \mathcal{Z}^{J}_{W}((v_{n+1}[m])_{k}\cdot\bfv_{n};B).
\end{align*}
\end{theorem}
\subsection{First reduction formula}
\label{SSfirst_red}

Suppose $L(0)v_{n+1}=\wt(v_{n+1})v_{n+1}$ and $J(0)v_{n+1}=\alpha v_{n+1}$,
$\alpha\in\C$.
Define $o_{\beta}(v_{n+1})=v_{n+1}(\wt(v_{n+1})-1+\beta)$ for $\beta\in\Z$.

\begin{lemma}[{\cite[Lem. 4.2]{BKT}}]
\label{lemRec1}
For all $\beta\in\Z$,
\[
  (1-\zeta^{-\alpha}q^{\beta}) \;T_{0}(o_{\beta}(v_{n+1})) \;
  \mathcal{Z}^{J}_{W}(\bfx_{n};B)
  =\sum_{k=1}^{n}\sum_{m\geq 0}
   \mathcal{Z}^{J}_{W}
   \;\left(\;\left(e^{z_{k}\beta}\frac{\beta^{m}}{m!}v_{n+1}[m]\cdot\right)_{\;k}
   \bfx_{n};B\right).
\]
\end{lemma}

\begin{corollary}[{\cite[Cor. 4.3]{BKT}}]
\label{corZeroRes}
Let $J(0)v_{n+1}=\alpha v_{n+1}$ and $\alpha z=\lambda\tau+\mu\in\Z\tau+\Z$.
Then
\[
  \sum_{k=1}^{n}\sum_{m\geq 0}
  \mathcal{Z}^{J}_{W}
  \;\left(\;\left(e^{z_{k}\lambda}\frac{\lambda^{m}}{m!}v_{n+1}[m]\cdot\right)_{\;k}
  \bfv_{n};B\right)=0.
\]
\end{corollary}

\begin{proposition}[{\cite[Prop. 4.4]{BKT}}]
\label{propZhured}
Assume $J(0)v_{n+1}=\alpha v_{n+1}$ and $\alpha z\notin\Z\tau+\Z$.
Then
\[
  \mathcal{Z}^{J}_{W}(\bfx_{n+1};B)
  =\sum_{k=1}^{n}\sum_{m\geq 0}
   \widetilde{P}_{m+1}\;\left(\frac{z_{n+1}-z_{k}}{2\pi i},\alpha z,\tau\right)
   \mathcal{Z}^{J}_{W}((v_{n+1}[m]\cdot)_{k}\bfx_{n};B).
\]
\end{proposition}

\begin{proposition}[{\cite[Prop. 4.5]{BKT}}]
\label{propZhured0}
Assume $J(0)v_{n+1}=\alpha v_{n+1}$ and $\alpha z=\lambda\tau+\mu\in\Z\tau+\Z$.
Then
\begin{align}
  \mathcal{Z}^{J}_{W}(\bfx_{n+1};B)
  &=e^{-z_{n+1}\lambda} \;\tr_{W}\;\bigl(v_{n+1}(\wt(v_{n+1})-1+\lambda) \;
    Y(\bfy_{n}) \;\zeta^{J(0)} \;q^{L(0)}\bigr)\notag\\
  &\quad+\sum_{k=1}^{n}\sum_{m\geq 0}
    P_{m+1,\lambda}\;\left(\frac{z_{n+1}-z_{k}}{2\pi i},\tau\right)
    \mathcal{Z}^{J}_{W}((v_{n+1}[m]\cdot)_{k}\bfx_{n};B),
  \label{eqZhuRed0}
\end{align}
with $P_{m+1,\lambda}$ as in \eqref{eqPellPm}.
\end{proposition}

\begin{proposition}[{\cite[Prop. 4.6]{BKT}}]
\label{propapnpt}
Assume $J(0)v_{n+1}=\alpha v_{n+1}$, $l\geq 1$, and
$\alpha z\notin\Z\tau+\Z$.
Then
\begin{eqnarray*}
  \mathcal{Z}^{J}_{W}(v_{n+1}[-l]\cdot x_{1},\bfx_{2,n};B)
  =\sum_{m\geq 0}(-1)^{m+1}\binom{m+l-1}{m}
    \widetilde{G}_{m+l}(\alpha z,\tau) \;
    \mathcal{Z}^{J}_{W}(v_{n+1}[m]\cdot x_{1},\bfx_{2,n};B)&&
\notag
\\
  +\sum_{k=2}^{n}\sum_{m\geq 0}(-1)^{l+1}\binom{m+l-1}{m}
    \widetilde{P}_{m+l}\;\left(\frac{z_{1}-z_{k}}{2\pi i},\alpha z,\tau\right)
    \mathcal{Z}^{J}_{W}(v_{n+1}[m]\cdot\bfx_{n};B).&&
\end{eqnarray*}
\end{proposition}

\begin{proposition}[{\cite[Prop. 4.7]{BKT}}]
\label{propapnpt0}
Assume $J(0)v_{n+1}=\alpha v_{n+1}$, $l\geq 1$, and
$\alpha z=\lambda\tau+\mu\in\Z\tau+\Z$.
Then
\begin{eqnarray*}
  \mathcal{Z}^{J}_{W}(v_{n+1}[-l]\cdot x_{1},\bfx_{2,n};B)
  =(-1)^{l+1}\frac{\lambda^{l-1}}{(l-1)!} \;
    \tr_{W}\;\bigl(v_{n+1}(\lambda+\wt(v_{n+1})-1) \;Y(\bfy_{n}) \;
    \zeta^{J(0)} \;q^{L(0)}\bigr)&&
\\
  +\sum_{m\geq 0}(-1)^{m+1}\binom{m+l-1}{m}
    E_{m+l,\lambda}(\tau) \;
    \mathcal{Z}^{J}_{W}(v_{n+1}[m]\cdot x_{1},\bfx_{2,n};B)&&
\\
  +\sum_{k=2}^{n}\sum_{m\geq 0}(-1)^{l+1}\binom{m+l-1}{m}
    P_{m+l,\lambda}\;\left(\frac{x_{1}-x_{k}}{2\pi i},\tau\right)
    \mathcal{Z}^{J}_{W}(v_{n+1}[m]\cdot\bfx_{n};B),&& 
\end{eqnarray*}
for $E_{k,\lambda}$ given in \eqref{eqGkl}.
\end{proposition}

\begin{remark}
For $\alpha=0$ one has $\lambda=\mu=0$, and
Propositions \ref{propZhured0} and \ref{propapnpt0} reduce to the standard
results of \cite{Zhu,MTZ} with $a(\lambda+\wt(a)-1)=o(a)$.
\end{remark}

\section{Vertex operator algebras}
\label{durdom}
\subsection{Vertex operator superalgebras}
\label{vertex}

We recall the axioms from \cite{B,FHL,FLM,K,MN}.
A  {vertex operator superalgebra} (VOSA) is a tuple
$(V,Y,\vac,\omega)$ where $V=V_{\bar{0}}\oplus V_{\bar{1}}=\bigoplus_{r\geq r_{0}}V_{r}$
($r_{0}\in\C$, $V_{r}=V_{\bar{0},r}\oplus V_{\bar{1},r}$),
$\vac\in V_{\bar{0},0}$ is the vacuum, $\omega\in V_{\bar{0},2}$ the Virasoro
vector, and
\[
  Y:V\longrightarrow(\End V)[[z,z^{-1}]],
  \quad
  Y(a)=\sum_{n\in\Z}a(n)z^{-n-1},
\]
satisfies: $a(n)\vac=\delta_{n,-1}a$ for $n\geq -1$;
$a(n)V_{\alpha}\subset V_{\alpha+p(a)}$;
and for all $x_{i}=(v_{i},z_{i})$,
\begin{equation}
\label{eqlocality}
  (z_{1}-z_{2})^{N}[Y(x_{1}),Y(x_{2})]=0\quad\text{for }N\gg 0
\end{equation}
(graded commutator).
The Virasoro operators $L(n)$ satisfy
\[
  [L(m),L(n)]=(m-n)L(m+n)+\tfrac{c}{12}(m^{3}-m)\delta_{m,-n}\Id_{V},
  \quad c\in\C,
\]
with $L(-1)$ the translation operator: $Y(L(-1)a)=\tfrac{d}{dz}Y(a,z)$;
$L(0)a=\wt(a)a$; $V_{r}=\{a\in V:\wt(a)=r\}$.
The standard commutator formula is
\begin{equation}
\label{eqcomm}
  [a(m),Y(x)]=\sum_{j\geq 0}\binom{m}{j}Y(a(j)\cdot x)z_{1}^{m-j}.
\end{equation}
For $a$ of weight $\wt(a)\in\Z$ the  {generalised zero mode} is
\[
  o_{\lambda}(a)=a(\wt(a)-1+\lambda),\quad\lambda\in\C,
\]
extended by linearity.
\subsection{Square-bracket formalism and shifted Virasoro vector}
\label{squa}

Following \cite{Zhu,DMs}, define
\[
  Y[x]=Y(e^{zL(0)}v,e^{z}-1)=\sum_{n\in\Z}v[n]z^{-n-1}.
\]
For $v$ of weight $\wt(v)$ \cite[Lem. 4.3.1]{Zhu},
\begin{equation}
\label{eqsq_bracket}
  \sum_{j\geq 0}\binom{k+\wt(v)-1}{j}v(j)=\sum_{m\geq 0}\frac{k^{m}}{m!}v[m].
\end{equation}
The square-bracket structure gives an isomorphic VOSA with Virasoro vector
$\widetilde{\omega}=\omega-\frac{c}{24}\vac$.

The  {shifted Virasoro vector} \cite{DMs} is
$\omega_{h}=\omega+h(-2)\vac$,
where $h=-\frac{\lambda}{\alpha}J$ for fixed $\lambda\in\Z$.
The shifted grading operator is $L_{h}(0)=L(0)+\frac{\lambda}{\alpha}J(0)$,
and the corresponding square-bracket vertex operators are
$Y[x]_{h}=Y(e^{zL_{h}(0)}v,e^{z}-1)=\sum_{n\in\Z}v[n]_{h}z^{-n-1}$,
satisfying
$Y[a,z]_{h}=e^{z\lambda}Y[a,z]$
and
$a[n]_{h}=\sum_{m\geq 0}\frac{\lambda^{m}}{m!}a[n+m]$.
\section*{Acknowledgements}
The author is supported by the
Institute of Mathematics, Academy of Sciences of the Czech Republic
(RVO 67985840). 

\end{document}